\documentclass [a4paper,12pt]{amsart}
\usepackage{graphicx}

\usepackage{amssymb}

\newtheorem{thm}{Theorem}[section]
\newtheorem{cor}[thm]{Corollary}
\newtheorem{lem}[thm]{Lemma}
\newtheorem{prop}[thm]{Proposition}

\theoremstyle{definition}
\newtheorem{defn}[thm]{Definition}

\theoremstyle{definition}
\newtheorem{rem}[thm]{Remark}

\theoremstyle{definition}
\newtheorem{exam}[thm]{Example}

\def\R{\mathbb R}

\def\Z{\mathbb Z}

\def\sgn{\operatorname{sgn}}
\def\ct{\operatorname{\hbox{$c_t$}}}
\def\nt{\operatorname{\hbox{$n_t$}}}
\def\ot{\operatorname{\hbox{$o_t$}}}
\def\At{\operatorname{\hbox{$A_t$}}}

\begin{document}

\author{A. Montesinos Amilibia and J.J.~Nu\~no Ballesteros}

\title{The self-linking number of a closed curve in $\R^n$}

\begin{abstract} We introduce the self-linking number of a smooth
closed curve $\alpha:S^1\to \R^n$ with respect to a
$3$-dimensional vector bundle over the curve, provided that some
regularity conditions are satisfied. When $n=3$, this construction
gives the classical self-linking number of a closed embedded curve
with non-vanishing curvature \cite{P}. We also look at some
interesting particular cases, which correspond to the osculating
or the orthogonal vector bundle of the curve.
\end{abstract}

\date{}

\address{Departament de Geometria i Topologia,
Universitat de Val\`encia, Campus de Burjassot, 46100 Burjassot
SPAIN}

\email{montesin@uv.es\quad nuno@uv.es}

\thanks{Work partially supported by DGICYT Grant PB96--0785}

\maketitle

\section{Introduction}
\begin{figure}[h]\label{helice}
\centerline{\includegraphics[scale=0.7]{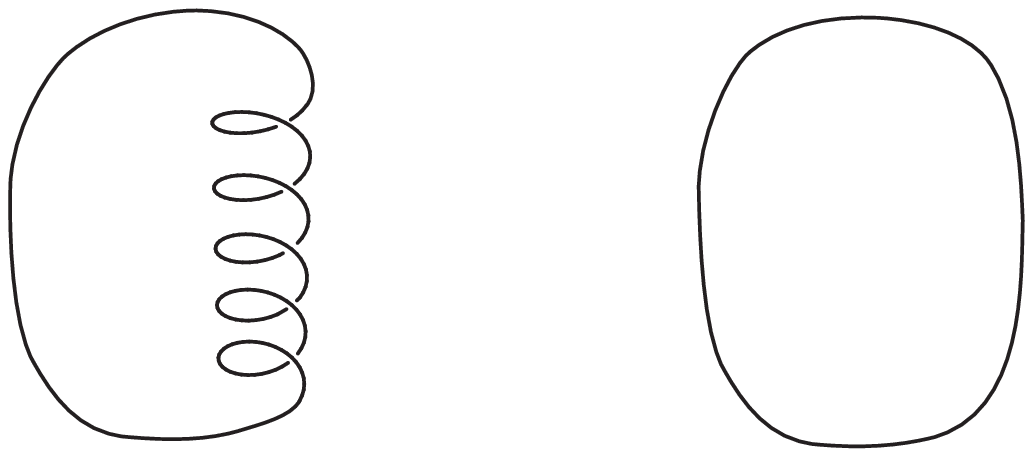}} \caption{}
\end{figure}

It is well known that two closed embedded curves
$\alpha,\beta:S^1\to\R^3$ are equivalent as knots if and only if
there is a continuous map $H:S^1\times[0,1]\to\R^3$ such that for
any $u\in[0,1]$, the curve $H_u:S^1\to\R^3$ given by
$H_u(t)=H(t,u)$ is an embedding and $H_0=\alpha, H_1=\beta$ (such
a map $H$ is said to be an \textit{isotopy} between
$\alpha,\beta$). For instance, if we look at the two curves shown
in Figure \ref{helice}, it follows that they are equivalent as
knots (in fact, they are equivalent to the trivial knot).

However, suppose that we construct these two curves so that they
are of class $C^3$ and have non-vanishing curvature at each point.
Then, it is not difficult to see that it is not possible to have
an isotopy $H$ such that for any $u\in[0,1]$, $H_u$ has the same
property (such a map will be called a \emph{non-degenerate
isotopy}). This is due to the fact that these two curves have
different self-linking number. This number was introduced by
C\u{a}lug\u{a}reanu \cite{C} and studied with more detail by Pohl
\cite{P}. It can be seen as the linking number between the given
curve and a curve obtained by slightly pushing the curve along the
principal normal. Moreover, it is possible to compute the
self-linking number by means of the following integral formula: $$
SL(\alpha)=\frac1{4\pi}\int_{S^1\times
S^1}\frac{\det(\alpha(s)-\alpha(t),\alpha'(s),\alpha'(t))}
{\|\alpha(s)-\alpha(t)\|^3}dt\wedge ds +\frac1{2\pi}\int_{S^1}
\tau dt,$$ where $\tau$ is the torsion of $\alpha$. Recently,
Gluck and Pan \cite{GP} have shown that there is a non-degenerate
isotopy between two embedded closed curves with non-vanishing
curvature in $\R^3$ if and only if they have the same knot type
and the same self-linking number. Thus, the self-linking number is
the key invariant if we want to do a curvature sensitive version
of the knot theory.

In this paper, we propose a generalization of this invariant for
the case of a closed smooth curve $\alpha:S^1\to\R^n$. In our
construction, we have to choose a $3$-dimensional vector bundle
over the curve, so that some regularity conditions hold between
the curve and the vector bundle. In the last part of the paper, we
analyze the osculating and the orthogonal self-linking number,
which correspond to the cases where the vector bundle is the
osculating or the orthogonal vector bundle of the curve,
respectively. These numbers can be also interpreted in terms of
intersection numbers of the curve with the orthogonal or the
osculating developable hypersurface of the curve, respectively. In
particular, it seems possible to relate them to some bitangency
properties of the curve \cite{N}.

A different approach in generalizing the self-linking number can
be found in \cite{W}, where it is considered a smooth map $f:M\to
\R^{2n+1}$ from a closed orientable smooth $n$-manifold $M$ into
$\R^{2n+1}$.

\section{The linking number of two curves with respect to a vector
bundle}\label{linking}

The linking number of two disjoint closed curves
$\alpha,\beta:S^1\to\R^3$ is a well known invariant, which is
defined as the degree of the map $e_1:S^1\times S^1\to S^2$ given
by $e_1(t,s)=(\beta(s)-\alpha(t))/\|\beta(s)-\alpha(t)\|$. In this
section, we will generalize this concept for two closed curves in
$\R^n$, by using a vector bundle over one of them.

Let $\nu:E\to S^1$ be a smooth 3-dimensional oriented vector
subbundle of the trivial vector bundle $S^1\times \R^n\to S^1$.
Given $t\in S^1$, we will denote the fiber by $\nu_t$, which is a
3-dimensional vector subspace of $\R^n$. Moreover, we will put
$\nt :\R^n\to \nu_t$ and $\ot :\R^n\to \nu_t^\bot$ for the
orthogonal projections, where $\nu_t^\bot$ is the orthogonal
subspace to $\nu_t$. Using these projections, we define the
\textit{covariant derivative} of a section $h:S^1\to\R^n$ of $\nu$
by $$(Dh)(t)=\nt (h'(t))=h'(t)-\ot (h'(t)),$$ and we will say that
$h$ is {\it parallel} if $Dh=0$.

Now, if we fix a parameterization of $S^1$ in the interval
$[0,\ell]$, we can solve the equations of parallel transport and
consider $\{p_1(t),p_2(t),p_3(t)\}$, a parallel oriented
orthonormal frame of $\nu_t$, for $t\in[0,\ell]$ (so that in
general $p_i(0)$ can be distinct from $p_i(\ell)$). This frame
allows us to define the linear map $\ct:\R^n\to\R^3$ by
$$\ct(x)=(\langle x,p_1(t)\rangle,\langle x,p_2(t)\rangle,\langle
x,p_3(t)\rangle),$$ so that the restriction of $\ct$ to $\nu_t$ is
an oriented isometry.

Given $h:S^1\to\R^n$ any smooth map, it will be useful to know
about the derivative of $\ct h(t)$. Let $\{u_i\}_{i=1}^3$ denote
the canonical basis of $\R^3$. Then,
\begin{align*}(\ct h(t))'&=\sum_{i=1}^3\langle
h(t),p_i(t)\rangle'u_i=\sum_{i=1}^3(\langle
h'(t),p_i(t)\rangle+\langle h(t),p_i'(t)\rangle)u_i\\ &=\ct
h'(t)+\sum_{i=1}^3\langle h(t),\ot p_i'(t)\rangle \ct
p_i(t)=\ct(h'(t)+\At h(t)),
\end{align*}
where $\At :\R^n\to\nu_t$ is the linear map given by $$\At
(x)=\sum_{i=1}^3\langle x,\ot p_i'(t)\rangle p_i(t).$$ It is not
difficult to see that $\At $ does not depend on the chosen
orthonormal frame $p_i(t)$, parallel or not.

\begin{defn} Let $\alpha,\beta:S^1\to\R^n$ be two smooth
closed curves in $\R^n$ and suppose that
$\beta(s)-\alpha(t)\notin\nu_t^\bot$ for any $(t,s)\in S^1\times
S^1$. We define the map $e_1:[0,\ell]\times S^1\to S^2$ by
$$e_1(t,s)=\frac{\ct(\beta(s)-\alpha(t))}{\|\ct(\beta(s)-\alpha(t))\|}.$$
Note that the imposed condition on the curves implies that
$\ct(\beta(s)-\alpha(t))\ne 0$ and thus, $e_1$ is well defined.
\end{defn}

\begin{lem}\label{forma} Let $\Omega_2$ be the standard volume form on $S^2$.
Then $e_1^*\Omega_2$ does not depend on the frame $p_i(t)$ and it
defines a closed smooth $2$-form on $S^1\times S^1$.
\end{lem}

\begin{proof} To abreviate, we denote
$\delta(t,s)=\beta(s)-\alpha(t)$. Then,
$$e_1^*\Omega_2=\frac{\det(\ct\delta(t,s),\partial_t
\ct\delta(t,s),\partial_s \ct\delta(t,s))}{\|
\ct\delta(t,s)\|^3}dt\wedge ds,$$ where $\partial_t$ and
$\partial_s$ denote the partial derivatives with respect to $t$
and $s$ respectively. But, according to the above computation,
$\partial_t \ct\delta(t,s)=\ct(-\alpha'(t)+\At \delta(t,s))$ and
$\partial_s \ct\delta(t,s)=\ct\beta'(s)$. Therefore,
\begin{align*}e_1^*\Omega_2&=\frac{\det(\ct\delta(t,s),
\ct(-\alpha'(t)+\At \delta(t,s)),\ct\beta'(s))}{\|
\ct\delta(t,s)\|^3}dt\wedge ds\\ &=\frac{\det(\nt \delta(t,s), \nt
(-\alpha'(t)+\At \delta(t,s)),\nt \beta'(s))}{\| \nt
\delta(t,s)\|^3}dt\wedge ds,\end{align*} where the last
determinant has to be considered with respect to any oriented
orthonormal frame of $\nu_t$.
\end{proof}

\begin{defn}\label{L} Let $\alpha,\beta:S^1\to\R^n$ be two smooth
closed curves in $\R^n$ and suppose that
$\beta(s)-\alpha(t)\notin\nu_t^\bot$ for any $(t,s)\in S^1\times
S^1$. We define its {\it linking number} with respect to $\nu$ as
$$L_\nu(\alpha,\beta)=\frac1{4\pi}\int_{S^1\times S^1}
e_1^*\Omega_2.$$
\end{defn}

It follows from this definition that if $\nu_t$ is constant, then
$L_\nu(\alpha,\beta)$ coincides with the classical linking number
of the projected curves, $L(c_t\circ \alpha, c_t\circ\beta)$. In
particular, when $n=3$, we have that
$L_\nu(\alpha,\beta)=L(\alpha,\beta)$, because then $\nu$ is the
trivial bundle.

\begin{lem} Let $\alpha,\beta,\nu$ be as in Definition~\ref{L}. Then,
there is $\tilde \beta:D^2\to \R^n$ extension of $\beta$, such
that $\tilde\beta(z)-\alpha(t)\notin\nu_t^\bot$ except for a
finite number of pairs $(t,z)\in S^1\times D^2$.
\end{lem}

\begin{proof} Let $\beta_1:D^2\to \R^n$ be an arbitrary extension
of $\beta$. Then the map $F:S^1\times D^2\times \R^n\to S^1\times
D^2\times \R^n$ given by $F(t,z,x)=(t,z,\beta_1(z)-\alpha(t)+x)$
is a diffeomorphism. In particular, it is transverse to the
submanifold $W=\{(t,z,x): x\in \nu_t^\bot\}$. By the
Transversality Theorem, it follows that for almost any $x\in\R^n$,
the map $F_x:S^1\times D^2\to S^1\times D^2\times \R^n$ given by
$F_x(t,z)=F(t,z,x)$ is also transverse to $W$. Since $W$ has
codimension 3, this implies that $F_x^{-1}(W)$ is finite.

To construct the required extension $\tilde \beta$, we piece
together $\beta_1$ near $S^1$ and $\beta_1+x$ on the interior as
follows. Let $\epsilon, \delta$ be such that
$0<\delta<\epsilon<1$. Let $g_{\epsilon,\delta}:D^2\to[0,1]$ be a
smooth function such that $g_{\epsilon,\delta}(z)=1$ if $\|z\|\le
\delta$ and $g_{\epsilon,\delta}(z)=0$ if $\epsilon\le \|z\|\le 1$
and let
$\beta_{\epsilon,\delta,x}(z)=\beta_1(z)+g_{\epsilon,\delta}(z)x$.
We claim that there are $\epsilon,\delta$ as above and $R>0$ such
that for any $x\in\R^n$ with $\|x\|<R$,
$$\beta_{\epsilon,\delta,x}(z)-\alpha(t)\in\nu_t^\bot\Longrightarrow
\beta_1(z)+x-\alpha(t)\in\nu_t^\bot.$$

Suppose that the claim is not true. Then, if for each $n>2$ we
consider $\epsilon=1-1/n$, $\delta=1-2/n$ and $R=1/n$, there are
$t_n\in S^1$, $z_n\in D^2$ and $x_n\in\R^n$ with $\|x_n\|<1/n$
such that
$$\beta_1(z_n)+g_{\epsilon,\delta}(z_n)x_n-\alpha(t_n)\in\nu_{t_n}^\bot,\text{
but } \beta_1(z_n)+x_n-\alpha(t_n)\notin\nu_{t_n}^\bot.$$ Thus
$g_{\epsilon,\delta}(z_n)\ne 1$ so that $\|z_n\|\ge 1-2/n$. By
taking subsequences if necessary, we can suppose that $t_n\to
t_0\in S^1$ and $z_n\to s_0\in S^1$. Thus, we arrive to
$\beta(s_0)-\alpha(t_0)\in\nu_{t_0}^\bot$, in contradiction with
the hypothesis. Now, we can choose
$\tilde\beta=\beta_{\epsilon,\delta,x}$, where $x$ is any one of
the points with $\|x\|<R$ for which $F_x^{-1}(W)$ is finite.
\end{proof}

\begin{prop} Let $\alpha,\beta,\nu$ be as in Definition~\ref{L}. Then,
$L_\nu(\alpha,\beta)\in\Z$.
\end{prop}

\begin{proof} Let $\tilde \beta$ be an extension of $\beta$ such
that $\tilde\beta(z)-\alpha(t)\notin\nu_t^\bot$ for any $(t,z)\in
S^1\times D^2\smallsetminus P$, being
$P=\{(t_1,z_1),\dots,(t_N,z_N)\}$. Then we can extend $e_1$ to
$\tilde e_1:[0,\ell]\times D^2\smallsetminus P\to S^2$ by putting
$$\tilde
e_1(t,z)=\frac{\ct(\tilde\beta(z)-\alpha(t))}{\|\ct(\tilde
\beta(z)-\alpha(t))\|}.$$ As in Lemma \ref{forma}, it follows that
$\tilde e_1^*\Omega_2$ defines a smooth $2$-form on $S^1\times
D^2\smallsetminus P$. Moreover, since $\Omega_2$ is closed on
$S^2$,  $d\tilde e_1^*\Omega_2=0$ and by Stokes Theorem,
$$0=\int_{S^1\times S^1} e_1^*\Omega_2+\sum_{i=1}^N\int_{\partial
B_i}\tilde e_1^*\Omega_2,$$ where $B_i$ denotes a small ball
centered at $(t_i,z_i)$ in the interior of $S^1\times D^2$ and
such that $B_i\cap B_j=\emptyset$ if $i\ne j$. In particular,
$$L_\nu(\alpha,\beta)=-\frac1{4 \pi}\sum_{i=1}^N\int_{\partial
B_i}\tilde e_1^*\Omega_2=-\sum_{i=1}^N\deg(\tilde e_1|_{\partial
B_i})\in\Z,$$ being $\deg(\tilde e_1|_{\partial B_i})$ the degree
of the map $\tilde e_1|_{\partial B_i}$.
\end{proof}

An immediate consequence of this, together with the fact that
$L_\nu(\alpha,\beta)$ depends continuously on $\alpha,\beta,\nu$
(when we consider the corresponding $C^\infty$ Whitney
topologies), is that $L_\nu(\alpha,\beta)$ is invariant under
homotopies of the curves and the vector bundle.

\begin{cor} Let $\alpha_u,\beta_u:S^1\to\R^n$ be
$1$-parameter families of curves and let $\nu_u:E_u\to S^1$ be a
$1$-parameter family of vector bundles, all of them depending
smoothly on the parameter $u\in[0,1]$ and such that
$\alpha_u,\beta_u,\nu_u$ satisfy the condition of
Definition~\ref{L}, for any $u\in[0,1]$. Then,
$L_{\nu_u}(\alpha_u,\beta_u)$ is constant on $u$.
\end{cor}

In the last part of this section, we give a characterization of
the linking number that will be used in the next section. Let
$\alpha,\beta,\nu$ be as in Definition~\ref{L} and suppose that
there is a vector field $\mu:S^1\to\R^n$ such that
$\mu(t)\notin\nu_t^\bot$, for any $t\in S^1$. Let
$\{f_i(t)\}_{i=4}^n$ be an orthonormal oriented frame of
$\nu_t^\bot$, that is, the basis
$(p_1(t),p_2(t),p_3(t),f_4(t),\dots,f_n(t)$ has the same
orientation as the canonical basis of $\R^n$. We can define the
map $\chi:S^1\times\R\times\R^{n-3}\to\R^n$ by
$$\chi(t,\lambda,x_4,\dots,x_n)=\alpha(t)+\lambda\mu(t)+\sum_{i=4}^n
x_i f_i(t).$$

\begin{prop}\label{indice} Suppose that $\beta$ meets the map $\chi$
transversely at a finite number of points and let
$$P_i=\beta(s_i)\in\alpha(t_i)+\lambda_i\mu(t_i)+\nu_{t_i}^\bot,\quad
i=1,\dots,N.$$ be those points. Then,
$$L_\nu(\alpha,\beta)=\frac12\sum_{i=1}^N\sgn(\lambda_i)
i(\beta,\chi;P_i),$$ where $i(\beta,\chi;P_i)$ denotes the
intersection number of $\beta$ and $\chi$ at $P_i$ and
$\sgn(\lambda_i)$ is the sign of $\lambda_i$.
\end{prop}

\begin{proof} Let $S^0=[0,\ell]\times S^1
\smallsetminus\{(t_1,s_1),\dots,(t_N,s_N)\}$. For any $(t,s)\in
S^0$, we have that $\ct\delta(t,s)\times \ct\mu(t)\ne 0$, where
$\delta(t,s)=\beta(s)-\alpha(t)$. Thus, we can define
$$e_3(t,s)=\frac{\ct\delta(t,s)\times
\ct\mu(t)}{\|\ct\delta(t,s)\times \ct\mu(t)\|}$$ and
$e_2(t,s)=e_3(t,s)\times e_1(t,s)$, so that $\{e_i(t,s)\}_{i=1}^3$
is a right-handed orthonormal frame of $\R^3$.

Now, we can consider the 1-forms on $S^0$ defined by
$\omega_{ij}=\left<de_i,e_j\right>$, for any $i,j=1,2,3$. Since
$\left<e_i,e_j\right>=\delta_{ij}$, by taking differentials we see
that $\omega_{ij}=-\omega_{ji}$. Moreover, we have that
\begin{align*} e_1^*\Omega_2&=\det(e_1,\partial_t e_1,\partial_s
e_1)dt\wedge ds\\&=(\omega_{12}(\partial_t)\omega_{13}(\partial_s)
- \omega_{12}(\partial_s)\omega_{13}(\partial_t))dt\wedge
ds\\&=\omega_{12}\wedge\omega_{13}.\end{align*} But from the fact
that $dde_i=0$, we deduce that $$0=\left<dde_i,e_j\right>=
d\omega_{ij}-\sum_{k=1}^3\omega_{ik}\wedge\omega_{kj}.$$ In
particular, $d\omega_{32}=\omega_{12}\wedge\omega_{13}
=e_1^*\Omega_2$ and it is not difficult to see that $\omega_{32}$
defines a 1-form on $S^1\times
S^1\smallsetminus\{(t_1,s_1),\dots,(t_N,s_N)\}$. This gives, by
Stokes Theorem, that $$L_\nu(\alpha,\beta)=\frac1{4\pi}
\int_{S^1\times S^1} e_1^* \Omega_2= \frac1{4\pi}\lim_{\epsilon\to
0}\sum_{i=1}^N\int_{\partial D_\epsilon(t_i,s_i)}\omega_{32},$$
where $D_\epsilon(t_i,s_i)$ denotes the disk centered at
$(t_i,s_i)$ of radius $\epsilon>0$ in $S^1\times S^1$. To conclude
the proof, we just have to show that for any $i=1,\dots, N$,
$$\frac1{2\pi}\lim_{\epsilon\to 0} \int_{\partial
D_\epsilon(t_i,s_i)}\omega_{32}=\sgn(\lambda_i)
i(\beta,\chi;P_i).$$

On one hand, if we put $m(t,s)=\ct\delta(t,s)\times \ct\mu(t)$, it
is easy to see that the left hand side is equal to $\pm 1$, in
accordance with the sign of $D=\det(\partial_t
m(t_i,s_i),\partial_s m(t_i,s_i),c_{t_i}\delta(t_i,s_i))$. If we
compute this, we get
\begin{align*}c_{t_i}\delta(t_i,s_i)&=\lambda_i c_{t_i}\mu(t_i)\\
\partial_s m(t_i,s_i)&=c_{t_i}\beta'(s_i)\times c_{t_i}\mu(t_i)\\
\partial_t m(t_i,s_i)&= c_{t_i}(-\alpha'(t_i)-\lambda_i\mu'(t_i)+A_{t_i}
(\delta(t_i,s_i)-\lambda_i\mu(t_i)))\times c_{t_i}\mu(t_i),
\end{align*}
and using the isometry between $\nu_t$ and $\R^3$,
$$D=\lambda_i\|n_{t_i}\mu(t_i)\|^2\det(n_{t_i}\beta'(s_i),
n_{t_i}(\alpha'(t_i)+\lambda_i\mu'(t_i)-A_{t_i}
(\delta(t_i,s_i)-\lambda_i\mu(t_i))),n_{t_i}\mu(t_i)).$$

On the other hand, if we suppose that
$\beta(s_i)=\chi(t_i,\lambda_i,x^i)$ for $x^i\in\R^{n-3}$, we have
that $i(\beta,\chi;P_i)$ is equal to $\pm1$ depending on the sign
of $$E=\det(\beta'(s_i),\partial_t
\chi(t_i,\lambda_i,x^i),\partial_\lambda
\chi(t_i,\lambda_i,x^i),\partial_{x_4}\chi(t_i,\lambda_i,x^i),
\dots,\partial_{x_n}\chi(t_i,\lambda_i,x^i)).$$ Now,
\begin{align*}\partial_t\chi(t_i,\lambda_i,x^i)&=
\alpha'(t_i)+\lambda_i\mu'(t_i)+\sum_{j=4}^nx_j^i f_j'(t_i),\\
\partial_\lambda\chi(t_i,\lambda_i,x^i)&=\mu(t_i),\\
\partial_{x_j}\chi(t_i,\lambda_i,x^i)&=f_j(t_i),
\end{align*}
and thus, $$E=\det(n_{t_i}\beta'(s_i),
n_{t_i}(\alpha'(t_i)+\lambda_i\mu'(t_i)+\sum_{j=4}^nx_j^i
f_j'(t_i)),n_{t_i}\mu(t_i)).$$ Finally, note that
$$c_{t_i}\sum_{j=4}^nx_j^i f_j'(t_i)=-c_{t_i}A_{t_i}
(\delta(t_i,s_i)-\lambda_i\mu(t_i)),$$ which implies the desired
result.
\end{proof}

\section{The self-linking number of a curve with respect to a vector bundle}

We shall define here the self-linking number of a smooth curve
$\alpha:S^1\to\R^n$ with respect to a vector bundle $\nu$, as the
linking number of $\alpha$ and $\tilde \alpha$, where $\tilde
\alpha:S^1\to\R^n$ is close enough to $\alpha$ and so that the
conditions of Definition~\ref{L} are satisfied. To ensure that
there exists such a curve $\tilde \alpha$, we need to assume that
$\alpha(s)-\alpha(t)\notin\nu_t^\bot$, for $s\ne t$. Moreover, we
also have to put some regularity conditions between the curve and
the fiber bundle on the diagonal $s=t$.

Throughout this section, we will suppose that $\alpha:S^1\to\R^n$
be a smooth closed curve in $\R^n$ and that $\nu$ is a smooth
3-dimensional oriented vector subbundle of the trivial vector
bundle, as in Section \ref{linking}.

\begin{lem}\label{lema-sl} Suppose that $\alpha,\nu$ satisfy the following conditions:
\begin{enumerate}
\item For any $s\ne t$, $\alpha(s)-\alpha(t)\notin\nu_t^\bot$.
\item There exists $1\le k\le n-2$ such that for any $t\in S^1$,
\begin{enumerate}
\item $\alpha'(t),\dots,\alpha^{(k+1)}(t)$ are
linearly independent;
\item $\alpha'(t),\dots,\alpha^{(k-1)}(t)\in\nu_t^\bot$;
\item $\langle\alpha^{(k)}(t)\rangle\oplus\langle\alpha^{(k+1)}(t)\rangle\oplus\nu_t^\bot$.
\end{enumerate}
\end{enumerate}
Then, there is $\delta_0>0$ such that
$\alpha(s)-\alpha_\delta(t)\notin\nu_t^\bot$, for any
$0<\delta<\delta_0$ and for any $(t,s)\in S^1\times S^1$, where
$\alpha_\delta(t)=\alpha(t)+\delta \alpha^{(k)}(t)$.
\end{lem}

\begin{proof} Suppose that this is not true. Then, for each $m\ge
1$, there are $\delta_m<1/m$ and pairs $(t_m,s_m)\in S^1\times
S^1$ such that
$\alpha(s_m)-\alpha(t_m)-\delta_m\alpha^{(k)}(t_m)\in\nu_{t_m}^\bot$.
By taking subsequences if necessary, we can suppose that $s_m\to
s_0\in S^1$ and $t_m\to t_0$. If $s_0\ne t_0$, we arrive to
$\alpha(s_0)-\alpha(t_0)\in\nu_{t_0}^\bot$, in contradiction with
condition 1. Otherwise, let $s_0=t_0$. If we denote by
$\{f_i(t)\}_{i=4}^n$ a frame for $\nu_t^\bot$, we have for any
$m\ge 1$, $$(\alpha(s_m)-\alpha(t_m))\wedge\alpha^{(k)}(t_m)\wedge
f_4(t_m)\wedge\dots\wedge f_n(t_m)=0.$$ Since
$$\alpha(s_m)=\alpha(t_m)+\sum_{j=1}^{k+1}\frac{\alpha^{(j)}(t_m)}{j!}(s_m-t_m)^j+O\left(
(s_m-t_m)^{k+2}\right),$$ we have after substitution and division
by $(s_m-t_m)^{k+1}$,
$$\alpha^{(k+1)}(t_m)\wedge\alpha^{(k)}(t_m)\wedge
f_4(t_m)\wedge\dots\wedge f_n(t_m)+O\left( s_m-t_m\right)=0.$$
This would imply that
$$\alpha^{(k+1)}(t_0)\wedge\alpha^{(k)}(t_0)\wedge
f_4(t_0)\wedge\dots\wedge f_n(t_0)=0,$$ in contradiction with
condition 2.(c).
\end{proof}

\begin{rem} When $n=3$, necessarily $k=1$ and $\nu_t^\bot=\{0\}$.
Thus, conditions 1 and 2 of Lema \ref{lema-sl} just say that
$\alpha$ is embedded and that $\alpha'(t),\alpha''(t)$ are
linearly independent, for any $t\in S^1$.
\end{rem}

\begin{defn}\label{sl} Suppose that $\alpha, \nu$ satisfy conditions 1 and 2 of Lemma
\ref{lema-sl} and consider $\alpha_\delta(t)=\alpha(t)+\delta
\alpha^{(k)}(t)$. The {\it self-linking number of $\alpha$ with
respect to $\nu$} is defined as $$SL_\nu(\alpha)=\lim_{\delta\to
0} L_\nu(\alpha_\delta,\alpha).$$

Note that since the linking number is invariant under homotopies,
Lemma \ref{lema-sl} ensures that $L_\nu(\alpha_\delta,\alpha)$
does not depend on $\delta$, if $\delta$ is small enough.
\end{defn}

\begin{figure}[h]
\centerline{\includegraphics[scale=1.0]{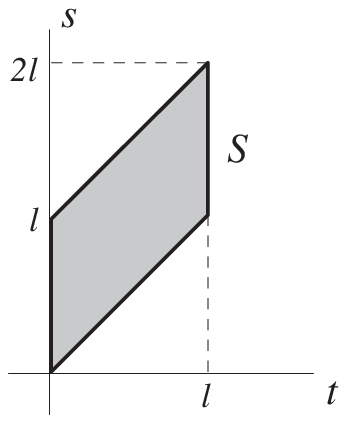}} \caption{}
\end{figure}

We would like now to obtain an integral expression for the
self-linking number analogous to the integral expression which
defines the linking number of two curves. The first step should be
to define the map $e_1$. Let $S$ be the following subset of $\R^2$
(see Figure 2) : $$S=\{(t,s)\in\R^2: 0\le t\le \ell,\ t\le s\le
t+\ell\}.$$ We define the map $e_1:S\to S^2$ as follows:
$$e_1(t,s)=\cases
\displaystyle{\frac{\ct(\alpha(s)-\alpha(t))}{\|\ct(\alpha(s)-\alpha(t))\|}},&
\text{ if $t<s<t+\ell$,}\\ \\
\displaystyle{\frac{\ct\alpha^{(k)}(t)}{\|\ct\alpha^{(k)}(t)\|}},&
\text{ if $s= t$,}\\
\\ (-1)^k\displaystyle{\frac{\ct\alpha^{(k)}(t)}{\|\ct\alpha^{(k)}(t)\|}},
& \text{ if $s= t+\ell$.}\endcases $$ Note that $e_1$ is well
defined when $\alpha$ satisfies conditions 1 and 2 of Lemma
\ref{lema-sl}. Moreover, by taking a Taylor expansion in a
neighbourhood of $(t,t)$ or $(t,t+\ell)$, it is easy to see that
$e_1$ is smooth.

\subsection{The case $k$ even}
If $k$ is even, the map $e_1$ can be considered as a map from
$[0,\ell]\times S^1$ to $S^2$. Moreover, $e_1^*\Omega_2$ defines a
closed 2-form on $S^1\times S^1$, which is the limit when
$\delta\to 0$ of the closed 2-form associated to the pair
$(\alpha_\delta,\alpha)$ in the definition of the linking number.
This gives the following result.

\begin{prop} Suppose that $\alpha, \nu$ satisfy
conditions $1$ and $2$ of Lemma
\ref{lema-sl} for $k$ even. Then,
$$SL_\nu(\alpha)=\frac1{4\pi}\int_{S^1\times S^1} e_1^*\Omega_2.$$
\end{prop}

\subsection{The case $k$ odd}
This case is more complicated. Let $S^0$ denote the open subset of
$S$ given by the pairs $(t,s)$ such that $t<s<t+\ell$ and $\ct
(\alpha(s)-\alpha(t))\times \ct\alpha^{(k)}(t)\ne0$. Then we can
complete $e_1$ on $S^0$ in order to get a frame of $\R^3$ as we
did in the proof of Proposition \ref{indice}. If $t<s<t+\ell$, we
define $$e_3(t,s)=\frac{\ct(\alpha(s)-\alpha(t))\times
\ct\alpha^{(k)}(t)} {\|\ct(\alpha(s)-\alpha(t))\times
\ct\alpha^{(k)}(t)\|}.$$ Moreover, it is possible to extend $e_3$
smoothly to the boundaries $s=t$ and $s=t+\ell$. In fact, by
taking a Taylor expansion in a neighbourhood of $(t,t)$ or
$(t,t+\ell)$, we get that
$$e_3(t,t)=e_3(t,t+\ell)=\frac{\ct\alpha^{(k+1)}(t)\times
\ct\alpha^{(k)}(t)} {\|\ct\alpha^{(k+1)}(t)\times
\ct\alpha^{(k)}(t)\|}. $$ We define $e_2$ in the obvious way,
$e_2(t,s)=e_3(t,s)\times e_1(t,s)$. Finally, we define the 1-forms
$\omega_{ij}=\left<de_i,e_j\right>$, for any $i,j=1,2,3$.

\begin{prop} Suppose that $\alpha, \nu$ satisfy
conditions $1$ and $2$ of Lemma \ref{lema-sl} for $k$ odd. Then,
$$SL_\nu(\alpha)=\frac1{4\pi}\int_{S^1\times S^1}
e_1^*\Omega_2-\frac1{2\pi}\int_{S^1} \phi,$$ where
$\phi(t)=\omega_{32}(t,t)$.
\end{prop}

\begin{proof} Let $\{f_i(t)\}_{i=4}^n$ be an orthonormal oriented frame of
$\nu_t^\bot$ and consider the map
$\chi:S^1\times\R\times\R^{n-3}\to\R^n$ given by
$$\chi(t,\lambda,x_4,\dots,x_n)=\alpha(t)+\lambda\alpha^{(k)}(t)+\sum_{j=4}^n
x_j f_j(t).$$ By the Transversality Theorem, we have that for a
residual subset of curves $\alpha$ and vector bundles $\nu$ with
the corresponding $C^\infty$ Whitney topologies, the curve
$\alpha$ meets the hypersurface $\chi$ transversely at a finite
number of points:
$$P_i=\alpha(s_i)=\alpha(t_i)+\lambda_i\alpha^{(k)}(t_i)+\sum_{j=4}^n
x_j^i f_j(t),\quad i=1,\dots,N,$$ with $s_i\ne t_i$ and
$\lambda_i\ne0$. Since $SL_\nu(\alpha)$, $\int_{S^1\times S^1}
e_1^*\Omega_2$ and $\int_{S^1} \phi$ depend continuously on
$\alpha$ and $\nu$, we can suppose that $\alpha$ and $\nu$ are
generic in the above sense.

In particular, for $\delta$ small enough, the same can be said if
we consider the intersection of $\alpha$ with $\chi_\delta$, where
$$\chi_\delta(t,\lambda,x_4,\dots,x_n)=\alpha_\delta(t)+\lambda\alpha^{(k)}(t)+\sum_{j=4}^n
x_j f_j(t).$$ Then, Proposition \ref{indice} gives that
$$L_\nu(\alpha_\delta,\alpha)=\frac12\sum_{i=1}^N\sgn(\lambda_i+\delta)
i(\alpha,\chi_\delta;P_i),$$ and taking limit when $\delta\to 0$,
$$SL_\nu(\alpha)=\frac12\sum_{i=1}^N\sgn(\lambda_i)
i(\alpha,\chi;P_i).$$

On the other hand, note that
$S^0=S\smallsetminus\{(t_1,s_1),\dots,(t_N,s_N)\}$. By using the
same argument as in the proof of Proposition \ref{indice},
$d\omega_{32}=\omega_{12}\wedge\omega_{13} =e_1^*\Omega_2$. If we
apply Stokes Theorem, $$\frac1{4\pi}\int_{S} e_1^* \Omega_2=
\frac1{4\pi}\int_{\partial S}\omega_{32}+\frac1{4\pi}
\lim_{\epsilon\to 0}\sum_{i=1}^N\int_{\partial
D_\epsilon(t_i,s_i)}\omega_{32},$$ where $D_\epsilon(t_i,s_i)$
denotes the disk centered at $(t_i,s_i)$ of radius $\epsilon>0$ in
$S$. Again we refer to the proof of Proposition \ref{indice} to
claim that $$\frac1{2\pi}\lim_{\epsilon\to 0}\int_{\partial
D_\epsilon(t_i,s_i)}\omega_{32}=\sgn(\lambda_i)
i(\alpha,\chi;P_i).$$ In particular,
$$SL_\nu(\alpha)=\frac1{4\pi}\int_{S}
e_1^*\Omega_2-\frac1{4\pi}\int_{\partial S}\omega_{32}.$$

To conclude the proof, we just have to compute the integral on
$\partial S$. We parameterize $\partial S$ by considering the
curves: $\gamma_1(u)=(0,u)$, $\gamma_2(u)=(u,u+\ell)$,
$\gamma_3(u)=(\ell,u+\ell)$ and $\gamma_4(u)=(u,u)$, for
$u\in[0,\ell]$. Then, we have that $$\int_{\partial S}\omega_{32}=
\int_0^\ell\omega_{32}(-\gamma'_1-\gamma'_2+\gamma'_3+\gamma'_4)
du. $$ Note that $\omega_{32}(0,u)=\omega_{32}(\ell,u+\ell)$ and
$\omega_{32}(u,u)=-\omega_{32}(u,u+\ell)$ for any $u\in[0,\ell]$.
This gives that $$\int_{\partial S}\omega_{32}=2\int_{S^1}\phi. $$
\end{proof}

\section{The orthogonal self-linking number}

We consider here the case that the vector bundle $\nu$ is equal to
the orthogonal vector bundle of the curve. That is, $\nu_t$ is the
$3$-plane orthogonal to the subspace generated by the $n-3$ first
derivatives of the curve.

\begin{defn}\label{SLO} Let $\alpha:S^1\to\R^n$ be a closed smooth curve in
$\R^n$ and suppose that:
\begin{enumerate}
\item For any $t\in S^1$, $\alpha'(t),\alpha''(t),\dots,\alpha^{(n-1)}(t)$
are linearly independent. In this way, at each point there is a
well defined Frenet frame $\{f_i(t)\}_{i=1}^n$ and also we have
the curvatures $\{\kappa_i(t)\}_{i=1}^{n-1}$. The {\it orthogonal
vector bundle} $\nu$ is defined so that $\nu_t$ is the $3$-plane
generated by $f_{n-2}(t),f_{n-1}(t),f_n(t)$.

\medskip
\item For any $s\ne t$ in $S^1$,
$\alpha(s)-\alpha(t)\notin\nu_t^\bot$. That is, the
$(n-3)$-osculating plane at $t$ does not meet the curve at any
other point.
\end{enumerate}
It follows that $\alpha,\nu$ satisfy conditions 1 and 2 of Lemma
\ref{lema-sl} for $k=n-2$. The self-linking number of $\alpha$
with respect to the orthogonal vector bundle will be called the
{\it orthogonal self-linking number} and will be denoted by
$SL^\bot(\alpha)$.
\end{defn}

With respect to this orthogonal vector bundle, we have that the
orthogonal projection $n_t:\R^n\to\nu_t$ is given by $$\nt(x)=
\left<x,f_{n-2}(t)\right>f_{n-2}(t)+\left<x,f_{n-1}(t)\right>
f_{n-1}(t)+\left<x,f_{n}(t)\right>f_{n}(t). $$ We also need to
know about the linear map $\At :\R^n\to\nu_t$. To simplify
computations, we will suppose that $\alpha$ is parameterized by
arc length. Then, $$\At (x)=\sum_{i=n-2}^n\langle x,\ot
f_i'(t)\rangle f_i(t)=-\langle x,f_{n-3}(t)\rangle
\kappa_{n-3}(t)f_{n-2}(t).$$ With this we can easily compute the
2-form $e_1^*\Omega_2$ used in the integral formula of the
self-linking number. But when $k=n-2$ is odd, we also need to
compute the $1$-form $\phi$ on $S^1$ given by
$\phi(t)=\omega_{32}(t,t)$. Note that
\begin{align*} e_1(t,t)&=\frac{\ct \alpha^{(k)}(t)}{\|\ct
\alpha^{(k)}(t)\|}=\ct f_{n-2}(t),\\ e_3(t,t)&=\frac{\ct
\alpha^{(k+1)}(t)\times \ct \alpha^{(k)}(t)}{\|\ct
\alpha^{(k+1)}(t)\times \ct \alpha^{(k)}(t)\|}=-\ct f_{n}(t),\\
e_2(t,t)&=e_3(t,t)\times e_1(t,t)=-\ct f_{n-1}(t),\\
\omega_{32}(t,t)&=\left<
{de_3(t,t)},e_2(t,t)\right>=-\kappa_{n-1}(t)dt.
\end{align*}
Thus, we have the following integral expression for the orthogonal
self-linking number.

\begin{cor}\label{integral-SLO} Let $\alpha:S^1\to\R^n$ be a closed smooth curve in
$\R^n$ satisfying conditions $1$ and $2$ of Definition \ref{SLO}.
Then
$$e_1^*\Omega_2=\frac{\left<\delta(t,s),f_{n-3}(t)\right>\kappa_{n-3}(t)
\det(\nt\delta(t,s),\nt\alpha'(s),f_{n-2}(t))}
{\|\nt\delta(t,s)\|^3}dt\wedge ds,$$ where
$\delta(t,s)=\alpha(s)-\alpha(t)$. Moreover, the orthogonal
self-linking number of $\alpha$ is equal to
$$SL^\bot(\alpha)=\cases\displaystyle{\frac1{4\pi}\int_{S^1\times
S^1}e_1^*\Omega_2},&\text{ when $n$ is even,}\\ \ \\
 \displaystyle{\frac1{4\pi}\int_{S^1\times
S^1}e_1^*\Omega_2 +\frac 1{2\pi}\int_{S^1}
\kappa_{n-1}(t)dt},&\text{ when $n$ is odd.}\endcases$$
\end{cor}

Given $\alpha:S^1\to\R^n$ a closed smooth curve in $\R^n$
satisfying conditions $1$ and $2$ of Definition \ref{SLO}, we can
consider the osculating developable hypersurface, which is the map
$\chi^\top:S^1\times\R^{n-2}\to\R^n$ defined by
$$\chi^\top(t,x_1,\dots,x_{n-2})=\alpha(t)+\sum_{i=1}^{n-2} x_i
f_i(t).$$ By condition 1, this is an immersion at those points
such that $x_{n-2}\ne 0$. Moreover, condition 2 implies that if
the curve meets this map at a point
$P=\alpha(s)=\chi^\top(t,x_1,\dots,x_{n-2})$ with $s\ne t$, then
necessarily $x_{n-2}\ne 0$. Note that the case $s=t$ would imply
that $x_1=\dots=x_{n-2}=0$.

\begin{cor} Suppose that $\alpha$ meets the map $\chi^\top$
transversely at a finite number of non-diagonal points and let
$(t_1,s_1),\dots,(t_N,s_N)$ be the pairs in $S^1\times S^1$
corresponding to these points. Then, the orthogonal self-linking
number of $\alpha$ is equal to
$$SL^\bot(\alpha)=-\frac12\sum_{i=1}^N
\sgn\langle\alpha'(s_i),f_n(t_i)\rangle.$$
\end{cor}

\section{The osculating self-linking number}

Here, we look at the self-linking number of a curve with respect
to its osculating vector bundle. That is, $\nu_t$ is the $3$-plane
generated by the $3$ first derivatives of the curve.

\begin{defn}\label{SLOR} Let $\alpha:S^1\to\R^n$ be a closed smooth curve in
$\R^n$ and suppose that:
\begin{enumerate}
\item For any $t\in S^1$, $\alpha'(t),\alpha''(t),\alpha'''(t)$
are linearly independent. We will denote by $\nu$ the osculating
vector bundle of $\alpha$, that is, $\nu_t=\langle
\alpha'(t),\alpha''(t),\alpha'''(t)\rangle$.

\medskip
\item For any $s\ne t$ in $S^1$,
$\alpha(s)-\alpha(t)\notin\nu_t^\bot$. That is, the
$(n-3)$-orthogonal plane at $t$ does not meet the curve at any
other point.
\end{enumerate}
In this case, $\alpha,\nu$ satisfy conditions 1 and 2 of Lemma
\ref{lema-sl} for $k=1$. Thus, we define the {\it osculating
self-linking number}, $SL^\top(\alpha)$, as the self-linking
number of $\alpha$ with respect to the osculating vector bundle.
\end{defn}

Now, we use $f_1(t),f_2(t),f_3(t)$ for the (partial) Frenet frame
of $\nu_t$ and $\kappa_1(t),\kappa_2(t)$ for the non-vanishing
curvatures. Then, $n_t:\R^n\to\nu_t$ is given by $$\nt(x)=
\left<x,f_{1}(t)\right>f_{1}(t)+\left<x,f_{2}(t)\right>
f_{2}(t)+\left<x,f_{3}(t)\right>f_{3}(t). $$ Again, we will
suppose for simplicity that $\alpha$ is parameterized by arc
length. Thus, $$\At (x)=\sum_{i=1}^3\langle x,\ot f_i'(t)\rangle
f_i(t)=k(t)\langle x,\ot \alpha^{(4)}(t)\rangle f_{3}(t),$$ where
$k(t)=1/\kappa_1(t)\kappa_2(t)$.

Finally,  we compute the $1$-form $\phi(t)=\omega_{32}(t,t)$:
\begin{align*} e_1(t,t)&=\frac{\ct \alpha'(t)}{\|\ct
\alpha'(t)\|}=\ct f_{1}(t),\\ e_3(t,t)&=\frac{\ct
\alpha''(t)\times \ct \alpha'(t)}{\|\ct \alpha''(t)\times \ct
\alpha'(t)\|}=-\ct f_{3}(t),\\ e_2(t,t)&=e_3(t,t)\times
e_1(t,t)=-\ct f_{2}(t),\\ \omega_{32}(t,t)&=\left<
{de_3(t,t)},e_2(t,t)\right>=-\kappa_{2}(t)dt.
\end{align*}
Thus, we have the following integral expression for the osculating
self-linking number.

\begin{cor}\label{integral-SLOR} Let $\alpha:S^1\to\R^n$ be a closed smooth curve in
$\R^n$ satisfying conditions $1$ and $2$ of Definition \ref{SLOR}.
Then $$e_1^*\Omega_2=\frac{
\det(\nt\delta(t,s),\nt\alpha'(s),\alpha'(t)-k(t)\langle
\delta(t,s),\ot \alpha^{(4)}(t) \rangle f_{3}(t))}
{\|\nt\delta(t,s)\|^3}dt\wedge ds,$$ where
$\delta(t,s)=\alpha(s)-\alpha(t)$. Moreover, the osculating
self-linking number of $\alpha$ is equal to $$SL^\top(\alpha)=
\frac1{4\pi}\int_{S^1\times S^1}e_1^*\Omega_2 +\frac
1{2\pi}\int_{S^1} \kappa_{2}(t)dt.$$
\end{cor}

Finally, we can compute the osculating self-linking number by
looking at the intersection of the curve with its orthogonal
developable. Let $\{f_j(t)\}_{j=4}^n$ be any orthonormal oriented
frame that trivializes $\nu_t^\bot$. We consider the orthogonal
developable hypersurface, which is the map
$\chi^\bot:S^1\times\R^{n-2}\to\R^n$ defined by
$$\chi^\bot(t,x_3,\dots,x_{n})=\alpha(t)+\sum_{i=3}^{n} x_i
f_i(t).$$ Since in this case $k=1$, we have by Definition \ref{sl}
that $SL^\top(\alpha)=\lim_{\delta\to 0} L_\nu(\alpha+\delta f_1
,\alpha)$. But it is not difficult to see that we obtain the same
number if we change $f_1$ by $f_2$ or $f_3$. Thus, we have the
following immediate consequence of Proposition \ref{indice}, for
$\mu(t)=f_3(t)$.

\begin{cor} Let $\alpha:S^1\to\R^n$ be a closed smooth curve in
$\R^n$ satisfying conditions $1$ and $2$ of Definition \ref{SLOR}.
Suppose that $\alpha$ meets the map $\chi^\bot$ transversely at a
finite number of non-diagonal points and let
$$P_i=\alpha(s_i)=\alpha(t_i)+\sum_{j=3}^{n} x_j^i f_j(t_i),\quad
i=1,\dots,N$$ be those points. Then, the osculating self-linking
number of $\alpha$ is equal to
$$SL^\top(\alpha)=\frac12\sum_{i=1}^N\sgn(x_3^i)
i(\alpha,\chi^\bot;P_i).$$
\end{cor}

\section{The examples}

In this last section, we will give some examples which show that
when $n>3$, the orthogonal and the osculating self-linking numbers
are not trivial and are independent. All the examples are in
$\R^4$ and the computations have been done with {\it Mathematica}
\cite{MN}. We compute the intersection of the curve with
$\chi^\top$ or $\chi^\bot$ and the corresponding indices.
Moreover, we also compute the integral value of $SL^\bot$ or
$SL^\top$ in order to ratify the results.

\begin{exam} Let $\alpha:S^1\to\R^4$ be the curve given by
$$\alpha(t)=\left( \cos (A + t) + {{\sin^2 (t)}},\cos (A +
2\,t),\cos (t), \frac {A \sin (3\,t)}{27} \right).$$ It follows
that for $A=1$ and $A=1.3$, the curve $\alpha$ satisfies
conditions 1 and 2 of Definition \ref{SLO} and Definition
\ref{SLOR}.

When $A=1$,  $\alpha$ meets $\chi^\bot$ transversely at four
points with indices $1,1,1,-1$. In fact, we compute numerically
the integral of Corollary \ref{integral-SLO} and obtain that
$SL^\top(\alpha)=1$.  If we look now at the intersection with
$\chi^\top$, there are just two points of transverse intersection,
both with index $1$. In this case, the integral formula of
Corollary \ref{integral-SLOR} gives $SL^\bot(\alpha)=1$.

When $A=1.3$, the intersection with $\chi^\bot$ gives again four
points with indices $1,1,1,-1$ and the numerical value of the
integral formula is $SL^\top(\alpha)=1$. However, although there
are two points of transverse intersection with $\chi^\top$, this
time the indices are $1,-1$ and the integral formula gives in this
case $SL^\bot(\alpha)=0$.
\end{exam}

\begin{exam} We consider now a different family of curves in
$\R^4$: $$\alpha(t)=\left( -\cos (A + t) + {\frac{A\,\sin
(2\,t)}{8}},
  {\frac{-{A^3}\,\cos (2\,t) }{8}} + \sin (A + t),
  {\frac{\sin (5\,t)}{125}},{\frac{{A^2}\,\sin (3\,t)}{27}}\right).$$
For $A=1.6$, $\alpha$ satisfies conditions 1 and 2 of Definition
\ref{SLO} and Definition \ref{SLOR}. The intersection with
$\chi^\bot$ is equal to six points, all of them having index $1$,
and the numerical computation of the integral formula gives
$SL^\top(\alpha)=3$. The intersection with $\chi^\top$ is also
equal to six points, but in this case two of them have index $1$
and the other four $-1$. The integral formula gives
$SL^\bot(\alpha)=-1$.
\end{exam}

\end{document}